\documentclass[12pt]{amsart}
\usepackage[dvips]{graphicx}
\usepackage{amsmath,graphics}
\usepackage{amsfonts,amssymb,hyperref}
\usepackage{xypic}
\theoremstyle{plain}
\newtheorem*{theorem*}{Theorem}
\newtheorem*{lemma*} {Lemma}
\newtheorem*{corollary*} {Corollary}
\newtheorem*{proposition*}{Proposition}
\newtheorem*{conjecture*}{Conjecture}
\newtheorem{theorem}{Theorem}[section]
\newtheorem{lemma}[theorem]{Lemma}
\newtheorem*{theorem1*}{Theorem 1}
\newtheorem*{theorem2*}{Theorem 2}
\newtheorem*{theorem3*}{Theorem 3}

\newtheorem{proposition}[theorem]{Proposition}
\newtheorem{conjecture}[theorem]{Conjecture}

\newtheorem{example}[theorem]{Example}

\theoremstyle{remark}
\newtheorem*{remark}{Remark}

\newtheorem{example*}{Example}
\newtheorem*{claim}{Claim}

\theoremstyle{definition}

\def\G{\Gamma}

\def\op{\operatorname}
\def\tipi{\ti{\pi}}
 \def\Q{\Bbb{Q}}  \def\Z{\Bbb{Z}} \def\R{\Bbb{R}} 
\def\N{\Bbb{N}}  \def\l{\lambda} \def\ll{\langle} \def\rr{\rangle}
 \def\a{\alpha}   \def\bp{\begin{pmatrix}}
\def\sm{\setminus} \def\ep{\end{pmatrix}} \def\bn{\begin{enumerate}} 
  \def\div{\mbox{div}} \def\en{\end{enumerate}}
\def\ba{\begin{array}} \def\ea{\end{array}}  
 \def\S{\Sigma}  \def\a{\alpha} \def\b{\beta} \def\ti{\tilde}

\def\be{\begin{equation}} \def\ee{\end{equation}} 
   
 \def\hom{\mbox{Hom}}  \def\gcd{\mbox{gcd}}

\def\co{\colon}
\def\zt{\Z[t^{\pm 1}]} \def\qt{\Q[t^{\pm 1}]}

\def\i{\iota}
\def\ol{\overline}

\begin{document}
\title{The Turaev and Thurston norms}
\author{Stefan Friedl}
\address{Fakult\"at f\"ur Mathematik\\ Universit\"at Regensburg\\   Germany}
\email{sfriedl@gmail.com}

\author{Daniel S. Silver}
\address{Department of Mathematics and Statistics\\ University of South Alabama}
\email{silver@southalabama.edu}

\author{Susan G. Williams}
\address{Department of Mathematics and Statistics\\ University of South Alabama}
\email{swilliam@southalabama.edu}

\date{\today}
\begin{abstract}
In 1986, W.\ Thurston introduced a (possibly degenerate) norm on the first  cohomology group of a $3$-manifold.
Inspired by this definition, Turaev introduced in 2002 a analogous norm on the first  cohomology group of a finite 2-complex. We show that if $N$ is the exterior of a link in a rational homology sphere, then the Thurston norm agrees with a suitable variation
of Turaev's norm defined on any 2-skeleton of $N$.
\end{abstract}
\maketitle

\section{Introduction}

In 1986, W.\ Thurston
 \cite{Th86} introduced a seminorm for $3$-manifolds $N$ with empty or toroidal boundary. 
It is a function $x_N\co H^1(N;\Q)\to \Q_{\geq 0}$ which
measures the complexity of surfaces that are dual to cohomology classes. 
We adopt the custom of referring to $x_N$ as the \textit{Thurston norm}. 
It plays a central role in $3$-manifold topology and we recall its definition
in Section \ref{section:thurstonnorm}, where we will also review several of its key properties.

Later, in 2002, V.\ Turaev \cite{Tu02} introduced an analogously defined seminorm for $2$-complexes.
For any finite $2$-complex $X$ with suitably defined boundary $\partial X$, V.\ Turaev defined 
$t_X\co H^1(X,\partial X;\Q)\to \Q_{\geq 0}$ using complexities of dual 1-complexes. Inspired by work of C.\ McMullen \cite{Mc02}, V.\ Turaev gave lower bounds for $t_X$ in terms of the multivariable  Alexander polynomial whenever the boundary of $X$ is empty.
The precise definition of $\partial X$
will be recalled in  Section \ref{section:turaevnorm}.
For the purpose of the introduction it suffices to know that if $N$ is a compact triangulated $3$-manifold,
then the $2$-skeleton $N^{(2)}$ is a finite $2$-complex with empty boundary.

A homotopy equivalence induces a canonical isomorphism of homology and cohomology groups which we use to identify the groups.  Examples given in \cite[p.~143]{Tu02} show that $t_X$ is \emph{not} invariant under homotopy.
We therefore introduce the following variation: For any finite $2$-complex $X$ with empty boundary, we define the \emph{Turaev complexity function} as follows.
If $\phi \in H^1(X;\Q)=\hom(\pi_1(X),\Q)$, then 
\[ \ol{t}_X(\phi):=\inf\left\{ t_{Y}(\phi\circ f) \,\left|\,\ba{l} Y\mbox{ is a finite $2$-complex with }\partial Y=\emptyset \mbox{ and } \\ f\co \pi_1(Y)\to \pi_1(X)\mbox{ is an isomorphism}\ea\right.\right\}.\]
Clearly $\ol{t}_X$  depends only on the fundamental group of $X$.
Since the minimum of two norms need not satisfy the triangle inequality, 
the Turaev function is not a seminorm, as we will see later in Proposition \ref{prop:notanorm}.

For any $3$-manifold $N$, we further define
\[\ol{t}_N(\phi):=\ol{t}_{N^{(2)}}(\phi),\]
where $N^{(2)}$ is the 2-skeleton of a triangulation of $N$. It is clear from the definition of $\ol{t}$ that $\ol{t}_N$ does not depend
on the choice of a triangulation.
\smallskip

Given a  3-manifold $N$, it is  natural to compare 
$x_N$ and $\ol{t}_{N}$ on $H^1(N;\Q)$.
In general, they do not agree. Indeed in  Section \ref{section:examplemore} we will see that there
exist many examples of closed $3$-manifolds $N$ and classes $\phi\in H^1(N;\Z)$  such that $\ol{t}_{N}(\phi)>x_N(\phi)$.
The underlying  reason is quite obvious: the Thurston norm is defined using complexities of surfaces,
whereas the Turaev function is defined using complexities of graphs. However, the complexity of a closed surface is lower by at least one than the complexity of any underlying $1$-skeleton.

It is therefore reasonable to restrict ourselves to the class of $3$-manifolds where Thurston norm-minimizing surfaces can always be chosen to have no closed component. 
In Lemma \ref{lem:linkinqhs3} we will see that if $N=\Sigma^3\sm \nu L$ is the exterior of a of a link $L$  in a rational homology sphere $\Sigma$, then $N$ has  this property. For simplicity of exposition we henceforth restrict ourselves to this type of 3-manifolds.

Using explicit and elementary constructions of $2$-complexes, we prove the following. 

\begin{theorem}\label{mainthm}
Let $N$ be the  exterior of a link in a rational homology sphere. Then
\[ \ol{t}_{N}(\phi)\leq x_N(\phi), \mbox{ for any }\phi \in H^1(N;\Q).\]
\end{theorem}

It is  natural to ask whether the extra freedom provided by working with 2-complexes instead of 3-manifolds allows us to get lower complexities.
Our main theorem says that this is not the case, at least if we restrict ourselves to irreducible link exteriors.
(Note that it follows from the definitions and the Sch\"onflies Theorem that the exterior of a link $L$ in $S^3$ is irreducible if and only if $L$ is non-split.)

\begin{theorem}\label{mainthm2}
Let $N$ be the  exterior of a link in a rational homology sphere. If $N$ is irreducible,
then
\[ \ol{t}_{N}(\phi)= x_N(\phi), \mbox{ for any }\phi \in H^1(N;\Q).\]
\end{theorem}

We will prove the inequality $\ol{t}_X(\phi) \geq x_{N}(\phi)$ by studying the Alexander norms of finite covers of $X$ and $N$,
and by applying the recent results of I.\ Agol \cite{Ag08,Ag13}, D. Wise \cite{Wi09,Wi12a,Wi12b}, P.\ Przytycki--D.\ Wise \cite{PW14,PW12}
and Y.\ Liu \cite{Li13}. We do not know of an elementary proof of Theorem \ref{mainthm2}. 

Theorem \ref{mainthm2} fits into  a long sequence of results showing that minimal-genus Seifert surfaces and Thurston norm-minimizing surfaces are `robust' in the sense that they `stay minimal' even if one relaxes some conditions. Examples of this phenomenon have been found by many authors, see for example \cite{Ga83,Ga87,Kr99,FV14,Na14,FSW13}.
\smallskip

The paper is organized as follows.
In Section \ref{section:definitions} we recall the definition of the Thurston and Turaev norms, and we introduce the
Turaev complexity function.
In Section \ref{section:lowerbounds} we discuss the Alexander norm for $3$-manifolds and $2$-complexes,
and we recall how they give lower bounds on the Thurston norm and Turaev complexity function, respectively.
In Section \ref{section:examplemore}, we first show that the Turaev complexity function
of the 2-skeleton can be greater than the corresponding Thurston norm.
We then show in Section \ref{section:alexfinitecovers} that the Thurston norm of  any irreducible $3$-manifold with non-trivial toroidal boundary
is detected by the Alexander norm of an appropriate finite cover.
Finally, in Section  \ref{section:proofmainthm} we put everything together to prove Theorem \ref{mainthm2}.

\subsection*{Conventions.}
All $3$-manifolds are compact, orientable and connected, and all $2$-complexes are connected,
unless it says specifically otherwise.

\subsection*{Acknowledgments.}
The first author gratefully acknowledges the support provided by the  SFB 1085 `Higher Invariants' at the Universit\"at Regensburg, funded by the Deutsche Forschungsgemeinschaft (DFG). 
The second and third authors thank the Simons Foundation for its support.
 We wish to thank the referee for carefully reading our paper and for pointing out a well-hidden mistake in an argument.

\section{The definition of the Thurston norm and the Turaev norm}\label{section:definitions}

\subsection{The Thurston norm and fibered classes}\label{section:thurstonnorm}

Let $N$ be a 3-manifold with empty or toroidal boundary.
The \emph{Thurston norm} of a class $\phi\in H^1(N;\Z)$ is defined as
 \[
x_N(\phi)=\min \{ \chi_-(\S)\, | \, \S \subset N \mbox{ properly embedded surface dual to }\phi\}.
\]
Here, $\chi_-(\S)$ is the complexity of a surface $\S$ with connected components $\S_1, \ldots ,\S_k$, given by
\[\chi_-(\S)=\sum_{i=1}^k \max\{-\chi(\S_i),0\}.\]
 Thurston
\cite{Th86} showed that $x_N$ defines a (possibly degenerate) norm on $H^1(N;\mathbb{Z})$.   Note  that any norm on $H^1(N;\Z)$
extends uniquely to a  norm on $H^1(N;\Q)$, which we denote by the same symbol.
\smallskip

We say that a class  $\phi\in H^1(N;\Q)$ is \emph{fibered} if 
there exists a fibration $p\co N\to S^1$ such that $\phi$ lies in the pull-back of $H^1(S^1;\Q)$ under $p$. By \cite{Ti70}, a  class $\phi \in H^1(N;\Q)$ is fibered if and only if it can be represented by a non-degenerate closed 1-form.

Thurston \cite{Th86}  showed the Thurston norm ball
\[ \{ \phi \in H^1(N;\Q)\,|\, x_N(\phi)\leq 1\}\]
is a polyhedron. This implies that if $C$ is a cone on a face of the polyhedron, then the restriction of $x_N$ to $C$ is a linear function.
To put differently, for any $\a,\b\in C$ and non-negative $r,s\in \Q_{\geq 0}$, the linear combination $r\a+s\b$ also lies in $C$, and $x_N(r \a + s \b) = r x_N(\a) + s x_N(\b)$.

Thurston \cite{Th86} also showed that any fibered class lies in the open cone on a top-dimensional face of the Thurston norm ball.
Furthermore, any other class in that open cone is also fibered.
Consequently, the set of fibered classes is the union of open cones on top-dimensional faces of the Thurston norm ball.
We will refer to these cones as the \emph{fibered cones of $N$.}
A class  $\phi\in H^1(N;\Q)$ in the closure of a fibered cone is \emph{quasi-fibered}.

\subsection{The Turaev norm and the  Turaev complexity function for 2--complexes}\label{section:turaevnorm}

As in \cite{Tu02}, a  \emph{finite 2--complex} is the underlying topological space of a finite connected 2-dimensional CW-complex such that each point has a neighborhood homeomorphic to the cone over a finite graph.
Examples of finite 2-complexes are given by compact surfaces (see \cite[p.~138]{Tu02}), 2-skeletons of finite simplicial spaces,
and the products of graphs with a closed interval.

The \emph{interior} of $X$, denoted $\op{Int} X$, is the set of points in $X$
 that have neighborhoods homeomorphic to $\R^2$. Finally the boundary  $\partial X$ of $X$ is the closure in $X$ of the set of all points of $X\sm \op{Int} X$
that have open neighborhoods in $X$ homeomorphic to $\R$ or to $\R\times \R_{\geq 0}$.
Note that $\partial X$ is a graph contained in the 1-skeleton of the CW-decomposition of X.
For example, if $X$ is a compact surface, then $\partial X$ is precisely the boundary of $X$ in the usual sense.

Following Turaev \cite{Tu02}, we say that a graph $\G$ in a finite $2$-complex is \emph{regular}
if $\G\subset X\sm \partial X$ and if there exists a closed neighborhood in $X\sm \partial X$ homeomorphic to $\G\times [-1,1]$ so that $\G=\G\times 0$.
A \emph{coorientation} for a regular graph $\G$  with components $\G_1,\dots,\G_k$ is
the choice of a component of $\G_i\times [-1,1]\sm \G_i,$ for each $i=1,\dots,k$.
A cooriented regular graph $\G\subset X$ canonically defines
an element $\phi_\G\in H^1(X,\partial X;\Z)$. Given any $\phi \in H^1(X,\partial X;\Z)$,
there exists a cooriented regular graph $\G$ with $\phi_\G=\phi$. (We refer to \cite{Tu02} for details.)
\smallskip

Let $X$ be a finite 2-complex with $\partial X=\emptyset$, and let $\phi \in H^1(X;\Z)$. The \emph{Turaev norm} of $\phi$ is 
\[ t_X(\phi):=\min \{ \chi_-(\G)\,|\, \G\subset X \emph{ cooriented regular graph with }
\phi_\G=\phi\},\]
where $\chi_-(\G)$ is the complexity of a graph $\G$ with connected components $\G_1,\dots,\G_k$, given by
\[ \chi_-(\G):=\sum_{i=1}^k \max\{-\chi(\G_i),0\}.\]
Turaev \cite{Tu02} showed that $t_X\co H^1(X;\Z)\to \Z_{\geq 0}$ is a (possibly degenerate) norm, and, as in the previous section,  $t_X$ extends to a norm
\[ t_X\co H^1(X;\Q)\to \Q_{\geq 0}.\]
In Theorem \ref{thm:disconnected} we will show that in general one has to allow disconnected graphs $\Gamma$ to minimize the Turaev norm.

As we already mentioned in the introduction, Turaev \cite[p.~143]{Tu02} showed that $t_X$ is in general not invariant under  homotopy equivalence.
(In fact Turaev showed that $t_X$ is not even invariant under simple homotopy.)
We therefore introduce a variation of the Turaev norm:
if $X$ is a finite $2$-complex with $\partial X=\emptyset$, then given $\phi \in H^1(X;\Q)=\hom(\pi_1(X),\Q)$ 
the \emph{Turaev complexity function} of $\phi$ is
\[ \ol{t}_X(\phi):=\inf\left\{ t_{\G}(\phi\circ f) \,\left|\,\ba{l} \G\mbox{ is a finite $2$-complex with }\partial \G=\emptyset \mbox{ and } \\ f\co \pi_1(\G)\to \pi_1(X)\mbox{ is an isomorphism}\ea\right.\right\}.\]
We make the following observations:
\bn
\item It is clear that  $\ol{t}_X$ is invariant under homotopy equivalence. In fact $\ol{t}_X$  depends only on the fundamental group of $X$.
\item Since $\ol{t}_X$ is the infimum of continuous homogeneous functions (i.e., functions with $f(\l x)=\l f(x)$ for $\l>0$),
$\ol{t}_X$ is upper semi-continuous and homogeneous.
\item The complexity function $\ol{t}_X$ is defined as the infimum of norms. Note that the minimum of two norms is in general no longer a norm. For example, the infimum of the two norms $a(x,y):=|x|$ and $b(x,y):=|y|$ on $\R^2$ is not a norm. We will see in Proposition \ref{prop:notanorm} that $\ol{t}_X(\phi)$  is, in general,  not a norm.
\item From the definition, it follows immediately that $\ol{t}_X(\phi) \leq t_X(\phi)$, for any $\phi \in H^1(X;\Q)$.
\item For any finite $2$-complex $X$, Turaev shows in \cite[Section~1.6]{Tu02} that $t_X$ is algorithmically computable. We do not know whether this is also the case for the Turaev complexity function $\ol{t}_X$. 
\en

\subsection{An inequality between the Thurston norm and the Turaev complexity function}
The goal of this section is to prove the following inequality between the Thurston norm and the Turaev complexity function.

\begin{proposition}\label{mainpropa}
Let $N$ be a 3-manifold and let $\phi \in H^1(N;\Z)$. 
If $\phi$ is dual to a properly embedded Thurston norm minimizing surface with $r$ closed components, then 
\[ \ol{t}_N(\phi) \leq x_{N}(\phi)+r.\]
\end{proposition}

\begin{proof}
Let  $\phi\in H^1(N;\Z)$ and let  $\S=\S_1\cup \dots\cup \S_s$ be 
 a surface  dual to $\phi$ of minimal complexity such that  $\S_1,\dots,\S_r$ are closed and $\S_{r+1},\dots,\S_s$ have nonempty boundary.

For $i=1,\dots,r$ we  pick an embedded graph $\Gamma_i\subset \Sigma_i$ with $\chi(\Gamma_i)=\chi(\Sigma_i)-1$ and such that $\pi_1(\Gamma_i)$ surjects onto $\pi_1(\Sigma_i)$. Furthermore, for $i=r+1,\dots,s$ we  pick an embedded graph $\Gamma_i\subset \Sigma_i$ with $\chi(\Gamma_i)=\chi(\Sigma_i)$ and such that $\pi_1(\Gamma_i)$ surjects onto $\pi_1(\Sigma_i)$.

Next we select pairwise disjoint product neighborhoods $\S_1\times [-1,1],\dots,\S_s\times [-1,1]$ such that the product orientations match the orientation of $N$.
We equip
\[ M:=N\,\,\sm \,\, \bigcup\limits_{i=1}^s \S_i\times (-1,1)\]
with a triangulation such that each $\G_i\times \{\pm 1\}$ is a subspace of $M^{(1)}$. 
Consider
\[ Y:=M^{(2)} \,\,\cup \,\,\bigcup\limits_{i=1}^s \G_i\times (-1,1).\]
It is straightforward
to see that $Y$ is a finite $2$-complex with $\partial Y=\emptyset$, and  the inclusion
map $Y\to N$ induces an isomorphism of fundamental groups. By slight abuse of notation we denote the restriction of $\phi$ to $Y$ again by $\phi$.

For $i=1,\dots,s,$ we  identify $\G_i$ with $\G_i\times 0$. It is clear that
$\G:=\G_1\cup \dots \cup \G_s$ is a regular graph on $Y$. Furthermore, with the obvious coorientation, we
 have $\phi_\G=\phi$.
It  follows that
\[ \ba{rcl} \ol{t}_N(\phi)\leq t_Y(\phi)\leq \chi_-(\G)&=&\sum\limits_{i=1}^r \max\{-\chi(\Gamma_i),0\}\,\,\,+\,\,\, \sum\limits_{i=r+1}^s \max\{-\chi(\Gamma_i),0\}\\
&\le&\sum\limits_{i=1}^r\max\{-\chi(\Sigma_i)+1,0\}\,\,\,+\,\,\, \sum\limits_{i=r+1}^s \max\{-\chi(\Sigma_i),0\}\\
&\le&\chi_-(\S)+r\\
&=& x_N(\phi)+r.\ea\]
\end{proof}

\noindent \textbf{Theorem \ref{mainthm}.}\emph{
Let $N$ be the  exterior of a link in a rational homology sphere.
Then for any $\phi\in H^1(N;\Q),$ we have
\[ \ol{t}_N(\phi) \leq x_{N}(\phi).\]}

\begin{proof}
Let $N$ be the  exterior of a link in a rational homology sphere.
We write $X=N^{(2)}$. Since $\ol{t}$ and $x_N$ are  homogeneous, it suffices to show that $\ol{t}_X(\phi) \leq x_{N}(\phi)$  for every $\phi\in H^1(N;\Z)$.
Assume that  $\phi\in H^1(N;\Z)$.
By  Lemma \ref{lem:linkinqhs3} (see Section \ref{section:examplemore}) there exists a Thurston norm-minimizing surface  dual to $\phi$ such that each component has nonempty boundary. The desired inequality  follows immediately from Proposition \ref{mainpropa}.
\end{proof}

\section{Lower bounds on the norms coming from Alexander polynomials}
\label{section:lowerbounds}

\subsection{The Alexander polynomial}

Let $X$ be a compact CW-complex, and let $\varphi\co$ $ H_1(X;\Z)\to H$ be a homomorphism onto a free abelian group.
We denote by $\ti{X}^\varphi$ the cover of $X$ corresponding to $\varphi\co \pi_1(X)\to H_1(X;\Z)\to H$.
The group $H$ is the deck transformation group of $\ti{X}^\varphi\to X$, and it acts on $H_1(\ti{X}^\varphi;\Z)$.
Thus we can view $H_1(\ti{X}^\varphi;\Z)$ as a $\Z[H]$-module.
Since $\Z[H]$ is a Noetherian ring, it  follows that $H_1(\ti{X}^\varphi;\Z)$ is a finitely presented $\Z[H]$-module.
This means that there exists an exact sequence
\[ \Z[H]^r\xrightarrow{A} \Z[H]^s\to H_1(\ti{X}^\varphi;\Z)\to 0.\]
After possibly adding columns of zeros, we can assume that $r\geq s$. Define the \emph{Alexander polynomial} of $(X,\varphi)$ to be
\[ \Delta_{X,\varphi}:=\mbox{gcd of all $s\times s$-minors of $A$.}\]
We refer to \cite{Fo54,Tu01,Hi12}  for the proof of the classical fact that $\Delta_{X,\varphi}$ is well-defined up to multiplication by a unit in $\Z[H]$, i.e., up to multiplication
by an element of the form $\epsilon h$, where $\epsilon \in \{-1,1\}$ and $h\in H$.

If $\varphi\co H_1(X;\Z)\to H:=H_1(X;\Z)/\mbox{torsion}$ is the  canonical projection,
 then we write $\Delta_X:=\Delta_{X,\varphi}$, and we refer to it as  the \emph{Alexander polynomial $\Delta_X$} of $X$.
 Furthermore, if $\phi\in H^1(X;\Z)=\hom(\pi_1(X),\Z)$,
 then we view the corresponding Alexander polynomial $\Delta_{X,\varphi}$ as an element in $\zt$
under the canonical identification of the group ring $\Z[\Z]$ with the Laurent polynomial ring $\zt$.

\subsection{The one-variable Alexander polynomials}
In this section we  relate the degrees of one-variable Alexander polynomials to the Thurston norm and to the Turaev complexity function.

In the following, given a non-zero polynomial $p(t)=\sum_{i=r}^s a_it^i$ with $a_r\ne 0$ and $a_s\ne 0$, we write
\[ \deg(p(t))=s-r.\]
Note that the degree of a non-zero one-variable Alexander polynomial is well-defined.

The following proposition is well known, see e.g., \cite{FKm06} for a proof.

\begin{proposition}\label{prop:th}
Let $N$ be a closed $3$-manifold and let $\phi \in H^1(N;\Z)$ be primitive. If $\Delta_{N,\phi}\ne 0$,  then
\[ x_N(\phi)\geq \deg(\Delta_{N,\phi})-2.\]
Furthermore, equality holds if $\phi$ is a fibered class and if $N\ne S^1\times S^2$.
\end{proposition}

We prove the following. 

\begin{proposition}\label{prop:tu}
Let $X$ be a finite $2$-complex with $\partial X=\emptyset$, and let be $\phi \in H^1(N;\Z)$ primitive. If $\Delta_{X,\phi}\ne 0$,  then
\[ \ol{t}_X(\phi)\geq \deg(\Delta_{X,\phi})-1.\]
\end{proposition}

\begin{proof}
Let $Y$ be a finite $2$-complex with $\partial Y=\emptyset$, and let $\psi \in H^1(Y;\Z)$ be primitive. If $\Delta_{Y,\phi}\ne 0$,  then it follows from the Claim 2 on page 152 of \cite{Tu02} that
\[ {t}_Y(\psi)\geq \deg(\Delta_{Y,\psi})-1.\]
The desired inequality 
\[ \ol{t}_X(\phi)\geq \deg(\Delta_{X,\phi})-1\]
is an immediate consequence of this fact and the observation that the Alexander polynomial  depends only on the fundamental group of $X$.
\end{proof}

\subsection{The Alexander norm}

Let $X$ be a compact connected CW-complex. We write $H:=H_1(X;\Z)/\mbox{torsion}$
and also $\Delta_X=\sum_{h\in H} a_hh$.
Let $\phi\in H^1(X;\Q)=\hom(\pi_1(X),\Q)=\hom(H,\Q)$. Following McMullen \cite{Mc02}, we define the \emph{Alexander norm} of $\phi$ by
\[ a_X(\phi):=\max\{ \phi(h)-\phi(g)\,|\, a_g\ne 0\mbox{ and } a_h\ne 0\}.\]
It is  straightforward to see that $a_X$ is indeed a norm on $H^1(X;\Q)$.
As in the proof of Proposition \ref{prop:tu}, we use that fact that the Alexander polynomial and thus the Alexander norm depend only on the fundamental group  of $X$.
More precisely, if $f\co Y\to X$ is a map of compact connected CW-complexes that induces an isomorphism of fundamental groups,
then
\[ f_*(\Delta_Y)=\Delta_X \in \Z[H_1(X;\Z)/\mbox{torsion}],\]
and thus, for any $f\in H^1(X;\Q)=\hom(\pi_1(X),\Q)$, we have
\be \label{equ:samea} a_Y(\phi\circ f^*)=a_X(\phi).\ee

We begin with the following theorem due to McMullen \cite{Mc02}.

\begin{theorem}\label{thm:mc02}\label{thm:mcm}
Let $N$ be a $3$-manifold  with empty or toroidal boundary and with $b_1(N)\geq 2$. Then
\[ a_N(\phi)\leq x_N(\phi)\mbox{ for any }\phi\in H^1(N;\Q).\]
Furthermore, equality holds for quasi-fibered classes.
\end{theorem}

\begin{proof}
Let $N$ be a $3$-manifold  with empty or toroidal boundary and with $b_1(N)\geq 2$.
McMullen \cite[Theorem~1.1]{Mc02} showed that
\[ a_N(\phi)\leq x_N(\phi)\mbox{ for any }\phi\in H^1(N;\Q)\]
and that equality holds for all integral fibered classes.
Since $a_N$ and $x_N$ are homogeneous, it  follows immediately that equality also holds
for all fibered classes and, in fact, for all  quasi-fibered classes.
\end{proof}

The following analogous theorem, which says that the Alexander norm also gives lower bounds on the Thurston norm and the Turaev complexity function,
is due to Turaev \cite{Tu02}.

\begin{theorem}\label{thm:tu02}
Let $X$ be a finite $2$-complex  with $b_1(X)\geq 2$ and such that
 $\partial X=\emptyset$. Then
\[ a_X(\phi)\leq \ol{t}_X(\phi)\leq t_X(\phi)\mbox{ for any }\phi\in H^1(X;\Q).\]
\end{theorem}

\begin{proof}
Let $Y$ be a finite $2$-complex  with $b_1(Y)\geq 2$ and such that
 $\partial Y=\emptyset$. Then by \cite[Theorem~3.1]{Tu02}, we have
\[ a_Y(\psi)\leq {t}_Y(\psi)\mbox{ for any }\psi\in H^1(Y;\Q).\]
The theorem now follows immediately from combining this result with the definition of $\ol{t}_X(\phi)$ and (\ref{equ:samea}).
\end{proof}

\section{Proofs}

\subsection{The Thurston norm and the Turaev complexity function for closed $3$-manifolds} \label{section:examplemore}
The combination  of  Propositions \ref{prop:th}, \ref{prop:tu} and  \ref{mainpropa} gives us the following theorem showing that the Thurston norm of a closed $3$-manifold
need not agree with  Turaev complexity function of its 2-skeleton.

\begin{theorem}\label{thm:norms}
Let $N\ne S^1\times S^2$ be a closed $3$-manifold and let $\phi \in H^1(N;\Z)$ be a primitive fibered class. Then
\[\ol{t}_N(\phi)=x_N(\phi)+1.\]
\end{theorem}

We also prove:

\begin{proposition}\label{prop:notanorm}
There exists a finite $2$-complex $X$ with $\partial X=\emptyset$ such that $\ol{t}_X$ does not satisfy the triangle inequality,
i.e., $\ol{t}_X$ is not a norm.
\end{proposition}

\begin{proof}
Let $N$ be a fibered $3$-manifold with $b_1(N)=2$. We write $X=N^{(2)}$ for some triangulation of $N$. As we mentioned in Section \ref{section:thurstonnorm},
by \cite{Th86}  there exists an open 2-dimensional cone $C\subset H^1(N;\Q)$
such that all classes in $C$ are fibered and such that $x_N$ is a linear function on $C$.

Given $\phi\in H^1(N;\Z)$ we denote by
\[ \op{div}(\phi):=\op{max}\{k\in \N\,|\, \mbox{there exists $\psi\in H^1(N;\Z)$ with $\phi=k\psi$}\}\]
the divisibility of $\phi$.
It follows from Theorem \ref{thm:norms} and the homogeneity of the Thurston norm and the Turaev complexity function that
\be \label{equ:tx} \ol{t}_X(\phi)=x_N(\phi)+\div(\phi)\mbox{ for any }\phi\in H^1(N;\Z)\cap C.\ee
We prove the following claim.

\begin{claim}
There exist $\a,\b\in C$ with $\div(\a)+\div(\b)<\div(\a+\b)$.
\end{claim}

Pick two primitive vectors $\phi,\psi\in C$ which are not colinear. 
Since $\phi$ and $\psi$ lie in the cone $C$, it follows that  any non-negative linear combination of $\phi$ and $\psi$ also lies in $C$. 

Select a coordinate system for $H^1(N;\Z)$, i.e., choose an identification of $H^1(N;\Z)$ with $\Z^2$. Since $\phi$ is primitive, we can assume that $\phi=(1,0)$. Since $\psi$ is also primitive, we know that $\psi=(x,y)$ for some coprime $x$ and $y$. Since $\phi$ and $\psi$ are not colinear, $y\ne 0$. 
Choose a prime $p>1+|y|$.  We  consider
$\a=(1,0)$ and $\b=(px+(p-1),py)$. Note that $p$ can not divide $px+p-1=p(x+1)-1$. It  follows that $\div(\b)=\gcd(px+(p-1),py)\leq |y|$. Evidently $\div(\a)=1$. Now 
\[ \div(\a+\b)\,=\,\div(px+p,py)\,=\,\gcd(px+p,py)\,\geq \,p\,> \, 1+|y|\,\geq \,\div(\a)+\div(\b).\]
This concludes the proof of the claim.

If we combine the claim and the linearity of $x_N$ on $C$ with equality (\ref{equ:tx}), then we obtain that
\[ \ba{rcl}\ol{t}_X(\a+\b)=x_N(\a+\b)+\div(\a+\b)&=&
x_N(\a)+x_N(\b)+\div(\a+\b)\\
&>&x_N(\a)+\div(\a)+x_N(\b)+\div(\b)\\
&=& \ol{t}_X(\a)+\ol{t}_X(\b).\ea\]
We have shown that $\ol{t}_X$ does not satisfy the triangle inequality.
\end{proof}

\subsection{The Alexander norm of finite covers of $3$-manifolds}\label{section:alexfinitecovers}

We begin with the following theorem. We state it in slightly greater generality than we actually need, since the result has independent interest.

\begin{theorem}\label{thm:agoletal}
Let $N\ne S^1\times D^2$ be an aspherical $3$-manifold  with empty or toroidal boundary. If $N$ is neither a Nil-manifold nor a Sol-manifold,  there exists a finite cover $p\co \ti{N}\to N$ such that $b_1(\ti{N})\geq 2$ and such that
\[ a_{\ti{N}}(p^* \phi)=x_{\ti{N}}(p^*\phi)\mbox{ for any $\phi \in H^1(N;\Q)$}.\]
\end{theorem}

The proof of the theorem will require the remainder of
Section  \ref{section:alexfinitecovers}.  
The theorem was proved for graph manifolds by Nagel \cite{Na14}. 
We will therefore restrict ourselves to the case of manifolds that are not (closed) graph manifolds. The main ingredient in our proof of Theorem \ref{thm:agoletal}
will be the following theorem, a consequence of the seminal work of Agol \cite{Ag08,Ag13}, Wise \cite{Wi09,Wi12a,Wi12b}, Przytycki-Wise \cite{PW14,PW12} and Liu \cite{Li13}. We summarize the main points of the proof for the convenience of the reader. 

\begin{theorem}\label{thm:vfib}
Let $N$ be an irreducible $3$-manifold  with empty or toroidal boundary that is not a closed graph manifold.
Then there exists a finite cover $p\co \ti{N}\to N$ such that,
for any $\phi \in H^1(N;\Q)$, the pull-back $p^*\phi$ is quasi-fibered.
\end{theorem}

\begin{proof}
Let $N$ be an irreducible $3$-manifold that is not a closed graph manifold.
It follows from the work of Agol \cite{Ag13}, Wise \cite{Wi09,Wi12a,Wi12b}, Przytycki-Wise \cite{PW14,PW12} and Liu \cite{Li13}
that $\pi_1(N)$ is virtually RFRS, i.e., $\pi_1(N)$ admits a finite index subgroup which is RFRS.
The precise definition of RFRS, references for which can be found in \cite{AFW15}, is not of concern to us.
What matters is that Agol \cite[Theorem~5.1]{Ag08} (see also \cite[Theorem~5.1]{FKt14}) showed that
if $\psi$ lies in $H^1(N;\Q)$ and if $N$ is an irreducible $3$-manifold such that $\pi_1(N)$ is virtually RFRS,
then there exists a finite cover $p\co \hat{N}\to N$ such that $p^*\psi$ lies in the closure of a fibered cone of $\hat{N}$.

By picking one class  in each cone of the Thurston norm ball of $N$ and iteratively applying Agol's theorem, one can easily show
that  there exists a finite cover $p\co\ti{N}\to N$ such that
for any $\phi \in H^1(N;\Q)$ the pull-back $p^*\phi$ lies in the closure of a fibered cone of $\ti{N}$. We refer to \cite[Corollary~5.2]{FV15}
for details.
\end{proof}

If $N$ is a graph manifold with nonempty boundary, then the conclusion of Theorem \ref{thm:vfib} also follows from facts that are more classical.
This argument is not used anywhere else in the paper, but since it is perhaps of independent interest we give a very quick sketch of the argument.

\begin{proof}[Proof of Theorem \ref{thm:vfib} if $N$ is a graph manifold]
Let $N$ be a graph manifold with boundary. It follows from Wang--Yu \cite[Theorem~0.1]{WY97}
and classical arguments, see e.g., \cite[Section~4.3.4.3]{AF13} and \cite{He87} that there exists a finite cover $\ti{N}$ of $N$
that is  fibered and such that if $\{N_v\}_{v\in V}$ denotes the set of  JSJ components of $\ti{N}$,
then each $N_v$ is of the form  $S^1\times \Sigma_v$ for some surface $\Sigma_v$.

For each $v\in V$ we write $t_v=S^1\times P_v$, where $P_v\in \Sigma_v$ is a point.
It follows from \cite[Theorem~4.2]{EN85} that a class $\phi\in H^1(\ti{N};\Q)$ is fibered if and only if $\phi(t_v)\ne 0$ for all $v\in V$.
Since $\ti{N}$ is fibered it now follows that all classes in $H^1(\ti{N};\Q)$ outside of finitely many hyperplanes are fibered.
Hence all classes in $H^1(\ti{N};\Q)$  are quasi-fibered.
\end{proof}

We can now move on to the proof of Theorem \ref{thm:agoletal}.
Note that arguments similar to the proof of Theorem \ref{thm:agoletal} were also used in \cite{FV15,FV14}.

\begin{proof}[Proof of Theorem \ref{thm:agoletal}]
Let $N\ne S^1\times D^2$ be an irreducible $3$-manifold  with empty or toroidal boundary that is not a closed graph manifold.
Since we assumed  that $N\ne S^1\times D^2$,  it now follows from Agol's Theorem \cite{Ag13}
and classical $3$-manifold topology that $N$ has a finite cover with $b_1$ at least two.  (We refer to \cite{AFW15} for details.)
We can therefore  assume that we already have $b_1(N)\geq 2$.

By Theorem \ref{thm:vfib} there exists a finite cover $p\co\ti{N}\to N$ such that
for any $\phi \in H^1(N;\Q)$, the pull-back $p^*\phi$ is quasi-fibered.
Note that Betti numbers never decrease by going to finite covers, i.e., we have $b_1(\ti{N})\geq b_1(N)\geq 2$.
It  follows from Theorem \ref{thm:mcm} that
\[ a_{\ti{N}}(p^* \phi)=x_{\ti{N}}(p^*\phi)\mbox{ for any $\phi \in H^1(N;\Q)$.}\]
This concludes the proof of the theorem.
\end{proof}

\subsection{Proof  of Theorem \ref{mainthm2}}\label{section:proofmainthm}

Before we turn to the proof of Theorem \ref{mainthm2}
we prove the following well-known lemma.

\begin{lemma} \label{lem:linkinqhs3}
If $N$ is  the exterior of a link in a rational homology sphere, then  any class $\phi\in H^1(N;\Z)$ is dual to a surface $\S$ of minimal complexity such that all components of $\Sigma$ have nonempty boundary. \end{lemma}


\begin{proof}
Let $N$ be the exterior of a link in a rational homology sphere.
It follows from a Mayer--Vietoris argument that the map $H_1(\partial N;\Q)\to H_1(N;\Q)$ is surjective. It follows from Poincar\'e duality and the Universal Coefficient Theorem that the boundary map $\partial \co H_2(N,\partial N;\Z)\to H_1(\partial N;\Z)$ has finite kernel. Since  $H_2(N,\partial N;\Z)\cong H^1(N;\Z)\cong \hom(H_1(N;\Z),\Z)$ is torsion-free it follows that the boundary map $\partial \co H_2(N,\partial N;\Z)\to H_1(\partial N;\Z)$ is in fact injective.
In particular this  implies  that closed surfaces represent the trivial homology class in $(N, \partial N)$.
Now let $\phi\in H^1(N;\Z)$, and let $\S$ be a properly embedded minimal-complexity surface dual to $\phi$. By the above observation, the closed components of $\S$ are null-homologous. It follows that the union of the components of $\S$ with non-trivial boundary represents the same homology as $\S$. Since removing components can never increase the complexity, we have shown that $\phi$ 
is dual to a surface $\S$ of minimal complexity such that all components of $\Sigma$ have nonempty boundary. 
\end{proof}

In the previous sections we collected all the tools that now allow us to finally complete the proof of Theorem \ref{mainthm2}.
\\

\noindent \textbf{Theorem \ref{mainthm2}.} 
\emph{ Let $N$ be the  exterior of a link in a rational homology sphere. If $N$ is irreducible, then for any $\phi\in H^1(N;\Q)$ we have}\[ \ol{t}_{N}(\phi) = x_{N}(\phi).\]

\begin{proof} It remains to prove that $ \ol{t}_{N}(\phi) \geq x_{N}(\phi)$.
Let $N$ be the  exterior of a link in a rational homology sphere.
Suppose that $N$ is irreducible. Let $\phi\in H^1(N;\Q)$.
It suffices to show that
if $Y$ is a   finite $2$-complex $Y$ with $\partial Y =\emptyset$ and if $f\co \pi_1(Y)\to \pi_1(N)$  is an isomorphism,
then
\[  t_{Y}(\phi\circ f)\geq x_N(\phi).\]
So let $Y$ and $f$ be as above. By a slight abuse of notation we denote $\phi\circ f\co \pi_1(Y)\to \Q$ by $\phi$ as well.

By Theorem \ref{thm:agoletal}  there exists a finite cover $p\co\ti{N}\to N$ such that $b_1(\ti{N})\geq 2$ and such that
\[ a_{\ti{N}}(p^* \phi)=x_{\ti{N}}(p^*\phi).\]
We  write $\pi=\pi_1(N)$ and $\ti{\pi}:=\pi_1(\ti{N})$,
and we denote by $p\co\ti{Y}\to Y$ the finite cover corresponding to $f^{-1}(\tipi)$. Note that $\ti{Y}$ is also a finite $2$-complex with $\partial \ti{Y}=\emptyset$. It  follows immediately from the definitions that
\[  x_{\ti{N}}(p^*\phi)\leq [\pi:\tipi]\cdot x_{N}(\phi)\mbox{ and } t_{\ti{Y}}(p^*\phi)\leq [\pi:\tipi]\cdot t_{Y}(\phi).\]
In fact, Gabai \cite[Corollary~6.13]{Ga83} showed that the above is an equality for the Thurston norm, i.e., we have the equality:
\[ x_{\ti{N}}(p^*\phi)=[\pi:\tipi]\cdot x_{N}(\phi).\]
Combining the above results with Theorem \ref{thm:tu02}, we   see that
\[ [\pi:\tipi]\cdot  {t}_{Y}(\phi) \geq {t}_{\ti{Y}}(p^*\phi)\geq a_{\ti{N}}(p^*\phi)=x_{\ti{N}}(p^*\phi)=[\pi:\tipi]\cdot x_{{N}}(\phi). \]
This concludes the proof the theorem.
\end{proof}

\subsection{Fundamental group complexity}\label{section:groupcomplexity}
Let $X$ be a finite $2$-complex with $\partial X = \emptyset$, and $\phi \in H^1(X; \Z)=
{\rm Hom}(\pi_1(X), \Z)$. In \cite{Tu02},
Turaev describes a method by which we can compute $t_X(\phi)$ using cocycles. We start by orienting edges (i.e., open $1$-cells) of $X$, and then select a $\Z$-valued cellular cocycle $k$ on $X$ representing $\phi$. We let 
\[ | k | = \sum_e (n_e/2-1)\vert k(e)\vert,\] 
where $e$ ranges over all edges in $X$, $k(e) \in \Z$ is the value of $k$ on $e$, and 
$n_e$ is the number of 2-cells adjacent to $e$, counted with multiplicity. (Note that
 $n_e\geq 2$ since $\partial X=\emptyset$.)  Turaev \cite[Section~1.6]{Tu02} proves that $t_X(\phi)$ is the minimum value of $\vert k\vert$ as $k$ ranges over all cellular cocycles representing $\phi$.
 
When the $0$-skeleton of $X$ consists of a single vertex, the 2-complex determines a group presentation $P$ for $\pi_1(X)$, and hence $|k|$ can be defined on the level of presentations: 

Given a finite presentation $P = \langle x_1, \ldots, x_m \mid r_1, \ldots, r_n\rangle $, following Turaev \cite{Tu02}, we denote by $\sharp (x_i)$  the  number of appearances 
of $x_i^{\pm 1}$ in the words $r_1, \ldots, r_n$. We say that $P$ is a \emph{good presentation}
if each $\sharp (x_i) \ge 2$. We are interested in good presentations, since it is  straightforward to see that  the canonical 2-complex corresponding to a good presentation has empty boundary. Also note that any finitely presented group admits a good presentation. Indeed, if $\# (x_i)=1$, then we can eliminate $x_i$ using a Tietze move. If $\# (x_i)=0$, then we can add a trivial relator $x_ix_i^{-1}$. 

Now let $P = \langle x_1, \ldots, x_m \mid r_1, \ldots, r_n\rangle $ be a good  presentation for a group $\pi$, and let $\phi$ be a homomorphism $\phi\colon \pi \to \Z$. We define  
\[t_P(\phi) = \sum_i (\sharp(x_i)/2-1)\,\vert \phi (x_i)\vert.\] 
Furthermore we define $\bar t_\pi(\phi)$ to be the minimum of $t_P(\phi)$ as $P$ ranges over all good
 presentations of $\pi$. We extend the definition in the usual way for rational cohomology classes $\phi \in H^1(X; \Q)$. 

\begin{lemma} \label{lem:gpcomplexity} 
Let $X$ be a finite $2$-complex with $\partial X= \emptyset$ and $\phi \in H^1(X; \Q)$.
We write $\pi=\pi_1(X)$.  Then 
\[\ol{t}_{X}(\phi) \le \bar t_\pi(\phi).\] 
\end{lemma} 

\begin{proof} Given a good
 presentation $P$ for $\pi$, we construct the
canonical finite 2-complex $Y$ with $\pi_1(Y) \cong \pi$.
 Let $k$ be the unique 1-cocycle representing $\phi$. A straightforward argument shows that  $\ol{t}_{X}(\phi) \le |k| = t_P(\phi)$, see also \cite[Section~1.8]{Tu02}.
  Since this is true for any good presentation of $\pi_1(X)$, we have $\ol{t}_{X}(\phi) \le \ol{t}_\pi(\phi)$. \end{proof} 

\begin{example}\label{wirtinger} Let $\pi$ the fundamental group of the exterior of a knot $K$ in the 3-sphere. Let $\phi$ be the abelianization homomorphism, mapping a meridian to $1$. If $P$ is a Wirtinger presentation corresponding to a diagram for $K$, then one  sees easily that $t_P(\phi)$ is the number of crossings of the diagram. 

It is usually possible to find presentations yielding a smaller value $t_P(\phi)$. Let $\Sigma$ be a Seifert surface for $K$ having minimal genus $g$. By splitting $\pi$ along $\pi_1(\Sigma)$, we obtain an HNN-decomposition for $\pi$ of the form \[\langle A, x \mid \mu(b) = x b x^{-1} \mbox{\ for all\ } 
b \in \pi_1(\Sigma) \rangle , \] where $A$ is the fundamental group of the knot exterior split along $\Sigma$, and
$\mu\colon \pi_1(\Sigma) \to A$ is injective. For such a presentation $P$, we have $t_P(\phi) = 2g-1$.  It follows by the next result that this value is the smallest possible; i.e., $\bar t_\pi(\phi) = 2g-1$.
 \end{example} 
 
\begin{theorem} \label{gpcomplexity} 
Let $N$ be the  exterior of a link in a rational homology sphere with group $\pi$.
If $N$ is irreducible, then for any $\phi \in H^1(N;\Q)$ such that $\Delta_{N, \phi} \ne 0$, we have
\[ \ol{t}_{N}(\phi) = \bar t_\pi(\phi)=  x_{N}(\phi).\]  \end{theorem}

\begin{remark} In \cite{Tu02}, Turaev gives several examples of knot groups and  presentations of minimal complexity. He states that it would be interesting to find other 
examples. Theorem \ref{gpcomplexity} shows how to construct presentations of minimal complexity for any knot in a rational homology sphere. \end{remark} 

\begin{proof} By Lemma \ref{lem:gpcomplexity} and Theorem \ref{mainthm2}, it suffices to prove that 
$\bar t_\pi (\phi) \le  x_{N}(\phi)$, for any $\phi \in H^1(N; \Q)$. 
By the homogeneity of the Turaev function and the Thurston norm we may assume that $\phi$ is an integral primitive cohomology class. 

Consider a Thurston norm-minimizing surface $\Sigma \subset N$ for $\phi$. Our assumption that $\Delta_{N, \phi}$ is not identically zero ensures that the first Betti number of ${\rm Ker}(\phi)$ is finite. 
 By a short argument in the beginning of the proof
 of \cite[Proposition~6.1]{Mc02}, the surface $\Sigma$ is connected. Its boundary is nonempty by Lemma \ref{lem:linkinqhs3}. Splitting $\pi$ along $\pi_1(\Sigma)$, as above, we obtain a presentation $P$ 
with complexity $2g-1$, where $g$ is the genus of $\Sigma$. Since $t_N(\phi) = 2g-1$, we are done. 
  \end{proof} 

We conclude this section with the following conjecture:

\begin{conjecture}
Let $X$ be a finite $2$-complex with $\partial X= \emptyset$. Then   
\[\ol{t}_{X}(\phi) = \bar t_{\pi_1(X)}(\phi) \mbox{ for  any $\phi \in H^1(X; \Q)$}.\] 
\end{conjecture}

Note that an affirmative answer to this question together with Theorems~\ref{mainthm} and \ref{mainthm2} would show that the conclusion of Theorem \ref{gpcomplexity} holds for any irreducible link complement $N$, without any assumptions on  $\phi$.

\section{Disconnected minimal dual graphs}
It is natural to ask whether one can always realize the Turaev norm of a primitive
 cohomology class by a  connected graph. 
In this final section of the paper we will see that this is not the case.
More precisely, we have the following theorem.

\begin{theorem}\label{thm:disconnected}
Given any $n$ there exists a 2-complex $X$ with $\partial X=\emptyset$ and a primitive class  $\phi \in H^1(\pi;\Z)$ such that for any 2-complex $Y$ with $\pi_1(Y)=\pi_1(X)$ and with $\partial Y=\emptyset$ the following holds:   any graph $\G$ in $Y$ that represents $\phi$ with $\ol{t}_X(\phi)=\chi_-(\G)$ has at least $n$ components.
\end{theorem} 

\begin{proof}
We consider the good presentation
\[ P=\ll a_1,\dots,a_n,x_1,\dots,x_n\,|\, [x_i,a_i], i=1,\dots,n\rr\]
and we denote by $X$ the corresponding 2-complex, which  is just  the join of $n$ tori $T_1,\dots,T_n$.
Clearly $\partial X=\emptyset$.

We write $\pi=\pi_1(X)$. The group $\pi$ is the free product of  $n$ free abelian groups $\ll a_i,x_i\,|\, [a_i,x_i]\rr$, $i=1,\dots,n$ of rank two. We consider  the epimorphism $\phi\co \pi\to \Z$ that is defined by $\phi(a_i)=0$, $i=1,\dots,n$ and $\phi(x_i)=1$, $i=1,\dots,n$.
 It is clear that on each torus $T_i$ there exists a circle, disjoint from the gluing point, such that the union of these circles is dual to $\phi$. We thus see that $\ol{t}_X(\phi)=0$. 

Now let $Y$ be a 2-complex  with $\pi_1(Y)=\pi$ and with $\partial Y=\emptyset$.
Let  $\Gamma$ be a graph on $Y$  which is dual to $\phi$ with $\chi_-(\Gamma)=0$.
We will show that  $\Gamma$ has at least $n$ components. 
Note that $\chi_-(\Gamma)=0$ implies that any component of $\Gamma$ is either a point or a circle. We denote by $m$ the number of components of $\Gamma$ that are circles. We will see that $m\geq n$. 

We start out with the following claim.

\begin{claim}
The module $H_1(Y;\qt)$ is isomorphic to $\qt^{n-1}\oplus \bigoplus_{i=1}^n \qt/(t-1)$.
\end{claim}

We first note that $H_1(Y;\qt)=H_1(X;\qt)$. A straightforward application of Fox
 calculus, see \cite{Fo53}, shows that  $H_1(X;\qt)\cong \qt^{n-1}\oplus \bigoplus_{i=1}^n \qt/(t-1)$. This concludes the proof of the claim.

Now we write $W=Y\sm \Gamma \times (-1,1)$. The usual Meyer--Vietoris sequence with $\qt$-coefficients corresponding to $Y=W\cup \Gamma\times [-1,1]$ gives rise to the exact sequence
\[ \dots \to H_1(\Gamma;\qt)\xrightarrow{\i_- -t\i_+} H_1(W;\qt)\to H_1(Y;\qt)\to H_0(\Gamma;\qt)\to \dots \]
Note that $\phi$ vanishes on $\Gamma$ and $W$. It follows that $H_*(\Gamma;\qt)$ and $H_*(W;\qt)$ are free $\qt$-modules. Furthermore, by the above discussion of $\Gamma$ we know that $H_1(\Gamma;\qt)\cong \qt^m$. It follows immediately from the above exact sequence and the classification of modules over PIDs that the torsion submodule of $H_1(Y;\qt)$ is generated by $m$ elements. 

On the other hand, we had just seen that the torsion submodule of  $H_1(Y;\qt)$ is isomorphic to $\oplus_{i=1}^n \qt/(t-1)$. It follows from the classification of modules over the PID $\qt$ that the minimal number of generators of the torsion submodule of  $H_1(Y;\qt)$ is $n$. Putting everything together we deduce that $m\geq n$. 
\end{proof}

\end{document}